\documentstyle[11pt]{article}
\begin {document}

\

\centerline {\Large \bf An intermediate value theorem} \centerline
{\Large \bf for sequences with terms in a finite set}

\

\

\centerline {\Large Mihai Caragiu}

\centerline {Department of Mathematics, Ohio Northern University}

\centerline {m-caragiu1@onu.edu}

\

\centerline {\Large Laurence D. Robinson}

\centerline {Department of Mathematics, Ohio Northern University}

\centerline {l-robinson.1@onu.edu}

\

\begin {abstract}
\noindent We prove an intermediate value theorem of an
arithmetical flavor, involving the consecutive averages $\{
\bar{x}_n\}_{n\geq 1}$ of sequences with terms in a given finite
set $\{ a_1,...,a_r\}$. For every such set we completely
characterize the numbers $\Pi $ ("intermediate values") with the
property that the consecutive averages $\{ \bar{x}_n\}$ of every
sequence $\{ x_n\}_{n\geq 1}$ with terms in $\{ a_1,...,a_r\}$
cannot increase from a value $\bar{x}_k<\Pi $ to a value
$\bar{x}_l>\Pi $ without taking the value $\bar{x}_s=\Pi$ for some
$s$ with $k<s<l$.
\\ \\ (2000) Mathematics Subject Classification: 11B99, 26D15
\end {abstract}

\

\

\newpage

\section {INTRODUCTION}

\

\noindent Let $r\geq 1$ be an integer and let $a_1,a_2,...,a_r$ be
real numbers with
$$a_1<a_2<...<a_r$$
Define
$$SEQ(a_1,...,a_r)$$
to be the set of all sequences $\{ x_n\}_{n\geq 1}$ such that
$x_n\in \{ a_1,a_2,...,a_r\}$ for all $n\geq 1$. For example
$SEQ(0,1)$ is the set of all binary sequences.

\

\noindent To each sequence $\{ x_n\}_{n\geq 1}$ we associate the
sequence of consecutive averages $\{ \bar{x}_n\}_{n\geq 1}$
defined by
$$\bar{x}_n=\frac {x_1+...+x_n}{n}$$
Clearly, if $\{ x_n\}_{n\geq 1}\in SEQ(a_1,...,a_r)$ then
$$a_1\leq \bar{x}_n\leq a_r$$
for all $n\geq 1$.

\

\noindent We are now in a position to define the sets which will
be studied in the current article.

\

\noindent {\bf DEFINITION.} For $a_1<a_2<...<a_r$ let us define
$$IV(a_1,...,a_r)$$
to be the set of all numbers $\Pi \in (a_1,a_r)$ with the
following "intermediate value property": {\it if $\{ x_n\}_{n\geq
1}\in SEQ(a_1,...,a_r)$ and if $\bar{x}_k<\Pi <\bar{x}_l$ for some
integers $k<l$ then there exists an integer $s$ with $k<s<l$ such
that $\bar{x}_s=\Pi $.}

\

\noindent A Putnam Exam problem [1] asks whether $\frac {4}{5}$ is
in $IV(0,1)$. Indeed the answer turns out to be affirmative. More
generally, we can make the following statement:

\

\noindent {\bf THEOREM 1.} $IV(0,1)=\left \{ \frac {k}{k+1}: k\geq
1 \right \} \subset (0,1)$.

\

\noindent In the present paper we will fully generalize the above
Theorem 1, providing a complete description of all "sets of
intermediate values" $IV(a_1,...,a_r)$. In particular we will
determine necessary and sufficient conditions under which
$IV(a_1,...,a_r)\neq \O$.

\

\noindent {\bf NOTE.} By definition, the numbers $\Pi \in
IV(a_1,...,a_r)$ are precisely those which cannot be "skipped" or
"jumped over" by {\it increasing} averages. In the last section we
will discuss the case of intermediate values which cannot be
skipped by {\it decreasing} averages. This being said, in the next
three sections, the term "skipped" will signify "skipped" by
averages going up (e.g. $\Pi = 0.7$ being skipped at the step
between the third and the fourth averages of the sequence
0,1,1,1,...).

\

\section {CASE OF BINARY SEQUENCES}

\

\noindent To prove THEOREM 1, we first show that
$$IV(0,1)\subseteq
\left \{ \frac {k}{k+1}: k\geq 1 \right \}.$$

\

\noindent Indeed, if
$$\Pi\in \left ( 0,\frac {1}{2} \right)\cup \left ( \frac {1}{2},\frac {2}{3}
\right) \cup \left ( \frac {2}{3},\frac {3}{4} \right) \cup \left
( \frac {3}{4},\frac {4}{5} \right)\cup ...$$ then the consecutive
averages of the sequence
$$0,1,1,1,1,...$$
will skip $\Pi $.

\

\noindent We now prove the reverse inclusion, namely,
$$\left \{ \frac {k}{k+1}: k\geq 1 \right \} \subseteq IV(0,1).$$
That is, we prove that if $\Pi =\frac {k}{k+1}$, $k\geq 1$ then
$\Pi $ cannot be skipped by a sequence of consecutive averages $\{
\bar{x}_n\}_{n\geq 1}$. We will proceed by contradiction. Assuming
$\{ \bar{x}_n\}_{n\geq 1}$ skips $\Pi =\frac {k}{k+1}$ it follows
that
$$\bar{x}_n<\frac {k}{k+1}<\bar{x}_{n+1}\leqno (1)$$
for some $n\geq 1$. First note that $x_{n+1}$ must be 1, because
if $x_{n+1}=0$ then the average $\bar{x}_{n+1}$ cannot be larger
than $\bar{x}_n$. Thus, if we denote $S=x_1+...+x_n$, (1) can be
rewritten as follows:
$$\frac {S}{n}<\frac {k}{k+1} < \frac {S+1}{n+1}.\leqno (2)$$
By cross-multiplying (2) is equivalent with the system of the
following two inequalities:
$$(k+1)S<nk,\leqno (3)$$
and
$$(n+1)k<(k+1)(S+1).\leqno (4)$$
From (3) and (4) it follows that
$$nk-1<(k+1)S<nk,$$
which is impossible, as all three terms are integers, and as there
can be no integer falling between consecutive integers. This
concludes the proof of THEOREM 1.

\

\noindent From THEOREM 1, a simple linearity argument leads us to
the following result.

\

\noindent {\bf THEOREM 2.} If $a<b$, then
$$IV(a,b)=\left \{ \frac {1}{k+1}a+\frac {k}{k+1}b: k\geq 1  \right
\}.$$

\

\noindent In the next section we will consider sequences with
terms in a set with three elements.

\

\section {CASE OF TERNARY SEQUENCES}

\

\noindent Let $0 < \mu < 1$. In  order to find $IV(0,\mu , 1)$ we
distinguish between the case of an irrational $\mu $ and the case
of a rational $\mu $. The easier case is the case of an irrational
$\mu $. Then there are no intermediate values for the consecutive
averages of sequences with terms in the set $\{ 0,\mu, 1\}$. In
other words:

\

\noindent {\bf THEOREM 3.} If $0< \mu < 1$ is irrational then
$IV(0,\mu , 1)=\O$.

\

\noindent {\bf PROOF.} We already know that every $\Pi \in (0,1)$
that is {\it not} of the form $\frac {k}{k+1}$ can be skipped by
the averages of some sequence in $SEQ(0,1)\subset SEQ(0,\mu , 1)$.
It will be enough to show that if $\mu $ is irrational then every
$\Pi \in (0,1)$ that is of the form $\frac {k}{k+1}$ can be
skipped by the averages of some sequence in $SEQ(0,\mu , 1)$.
Indeed, every $\Pi =\frac {k}{k+1} < \mu $ will be skipped by the
consecutive averages of the sequence
$$0,\mu, \mu, \mu,...$$ (the averages form an increasing sequence of irrationals
with limit $\mu $), while every $\Pi =\frac {k}{k+1}> \mu $ will
be skipped by the consecutive averages of the sequence
$$\mu ,1,1,1,...$$
(here, the averages form an increasing sequence of irrationals
with limit $1$). This concludes the proof of THEOREM 3.

\

\noindent Next we consider the case $0<\mu =\frac {p}{q}<1$ with
$p,q$ relatively prime positive integers. We will prove the
following result:

\

\noindent {\bf THEOREM 4.} If $0<\frac {p}{q}<1$ with $p,q$
relatively prime, then
$$IV(0,\frac {p}{q},1)= \left \{ 1-\frac {1}{qt}:t=1,2,3,... \right \}.$$

\

\noindent {\bf PROOF.} We know that every $\Pi \in (0,1)$ that is
{\it not} of the form $\frac {k}{k+1}$ can be skipped. The
question that remains is which numbers of the form $\Pi = \frac
{k}{k+1}$ can be skipped by the consecutive averages of some
sequence in $SEQ(0,\mu , 1)$.

\

\noindent First note that the sequence
$$x_1=\frac {p}{q},x_2=x_3=...=1$$
has consecutive averages of the form $$\frac {\frac
{p}{q}+l}{l+1}=\frac {p+ql}{q+ql},\leqno (5)$$ where
$l=0,1,2,...$. Then every $\Pi =\frac {k}{k+1}> \frac {p}{q}$ that
is not of the form $\frac {p+ql}{q+ql}$ can be skipped by the
sequence of increasing averages (5), so it cannot be in
$IV(0,\frac {p}{q},1)$. We now need to determine which fractions
of the form $\frac {p+ql}{q+ql}$ are also of the form $\frac
{k}{k+1}$. The equality
$$\frac {p+ql}{q+ql}=\frac {k}{k+1}$$
can be rewritten in the following equivalent form:
$$p(k+1)=q(k-l).\leqno(6)$$
From (6), keeping in mind that $p,q$ are relatively prime, we get
$$k+1=qt,$$
and
$$k-l=pt$$
for some integer $t$. In this case, $\frac {k}{k+1}$ is of the
form $\frac {qt-1}{qt}=1-\frac {1}{qt}$. As a consequence, if a
$\Pi = \frac {k}{k+1} > \frac {p}{q}$ is in the set of
intermediate values $IV(0,\frac {p}{q},1)$ then $q$ divides $k+1$.

\

\noindent We will now prove that there is no $\Pi $ in $IV(0,\frac
{p}{q},1)$ that is of the form $\frac {k}{k+1} $ and is less than
$ \frac {p}{q}$. To do this, we will show that if
$$\Pi =\frac {k}{k+1}<\frac {p}{q}$$
then there exists a sequence in $SEQ(0,\frac {p}{q},1)$ whose
averages skip $\Pi $. Indeed, let us consider the sequence
$$x_1=0,x_2=x_3=...=\frac {p}{q}. \leqno (7)$$
The averages of the sequence (7) form an increasing sequence
approaching $\frac {p}{q}$. If no $\bar {x}_k$ equals $\Pi $, then
$\Pi $ will be skipped. If $\bar {x}_n=\Pi $ for some $n$, then we
may consider the sequence
$$x_1=0,x_2=x_3=...=x_{n-1}=\frac {p}{q},x_n=1\leqno (8)$$
obtained by changing the $n$-th term of (7) into a one. Clearly,
$\Pi $ will be skipped at the transition between the $n-1$-th and
the $n$-th averages of the sequence (8).

\

\noindent At this point we know that every $\Pi \in (0,1)$ that is
{\it not} of the form $1-\frac {1}{qt}$ can be skipped by the
consecutive averages of some sequence in $SEQ(0,\frac {p}{q},1)$,
in other words,
$$IV(0,\frac {p}{q},1)\subseteq \left \{ 1-\frac {1}{qt}:t\geq 1 \right \}.\leqno (9)$$
The reverse inclusion
$$\left \{ 1-\frac {1}{qt}:t\geq 1 \right \}\subseteq IV(0,\frac {p}{q},1) .\leqno (10)$$
will be proved by contradiction. Assume that $\frac {qt-1}{qt}$
can be skipped by the consecutive averages of some sequence in
$SEQ(0,\frac {p}{q},1)$. Without loss of generality we may assume
that $\frac {qt-1}{qt}$ is in between the average
$$\bar{x}_n= \frac {x_1+...+x_n}{n}$$
with $u$ of the $x_1,...,x_n$ being zeros, $v$ being $\frac
{p}{q}$ and $w$ being ones ($u+v+w=n$) and the average
$$\bar{x}_{n+1}= \frac {x_1+...+x_n+1}{n+1}$$
with $u$ of the $x_1,...,x_n$ being zeros, $v$ being $\frac
{p}{q}$ and $w+1$ being ones ($x_{n+1}=1$):
$$\frac {v\frac {p}{q}+w}{n}<\frac {qt-1}{qt}<\frac {v\frac {p}{q}+w+1}{n+1}.\leqno (11)$$
Equivalently, (11) can be rewritten as follows:
$$\frac {pv+qw}{n}<\frac {qt-1}{t}<\frac {pv+qw+q}{n+1},$$
which is equivalent with the system consisting of the following
two inequalities:
$$pvt+qwt < nqt-n$$
and
$$nqt-n+qt-1 < pvt + qwt + qt$$
By using $n = u+v+w$ and after simplifying, the above two
inequalities will be equivalent to the following system:
$$pvt < qut + qvt -n\leqno (12)$$
and
$$qut + qvt-n-1 < pvt.\leqno (13)$$
From (12) and (13) it follows that
$$pvt<qut+qvt-n<pvt+1,$$
which is again a contradiction (as there can be no integer falling
between consecutive integers). This concludes the proof of the
reverse inclusion (10). From (9) and (10), THEOREM 4 follows.

\

\noindent A straightforward linearity argument based on the
previous two theorems leads us to the following intermediate value
theorem characterizing all sets $IV(a,b,c)$.

\

\noindent {\bf THEOREM 5.} Let $a<b<c$ and let $\mu : =\frac
{b-a}{c-a}$. If $\mu $ is irrational then
$$IV(a,b,c)=\O.$$
If $\mu =\frac {p}{q}$ with $p,q$ relatively prime positive
integers, then
$$IV(a,b,c)=\left \{ \left ( 1-\frac {1}{qt} \right )c+\frac {1}{qt}a:t=1,2,3,... \right
\}$$

\

\section {THE GENERAL INTERMEDIATE VALUE \\ THEOREM}

\

\noindent We will now completely characterize the intermediate
value sets of the form $IV(0,\mu_1,...,\mu _r, 1)$ where
$0<\mu_1<...<\mu _r<1$. First, note we can immediately obtain the
following result.

\

\noindent {\bf THEOREM 6.} If $\mu _i$ is irrational for some
$i=1,2,...,r$, then
$$IV(0,\mu_1,...,\mu _r, 1)=\O.$$

\

\noindent {\bf PROOF.} Follows from THEOREM 3, since if $\mu _i$
is rational then every $\Pi \in (0,1)$ can be skipped by the
averages of some sequence in $$SEQ(0,\mu_i,1)\subset
SEQ(0,\mu_1,...,\mu_r,1).$$

\

\

\noindent Now let us assume that all $\mu _i$'s are rational:
$$\mu _i=\frac {p_i}{q_i},\ \ i=1,..,r,$$
with $\gcd (p_i,q_i)=1$ for all $i=1,...,r$.

\

\noindent Let $M$ be the least common multiple of the denominators
of the reduced fractions $\frac {p_i}{q_i}$, $i=1,...,r$. We will
prove that the following result holds true.

\

\noindent {\bf THEOREM 7.} With the above notations, we have
$$IV(0,\mu_1,...,\mu _r, 1)=\left \{ 1-\frac {1}{Mt}:t=1,2,3,... \right
\} \leqno (14)$$

\

\noindent {\bf PROOF.} Let $i\in \{ 1,2,...,r\}$. From THEOREM 4
it follows that for $A\geq 2$, the element
$$\Pi = 1-\frac {1}{A}$$
will be skipped by the averages of some sequence in $$SEQ_(0,\frac
{p_i}{q_i},1)\subset SEQ(0,\mu_1,...,\mu_r,1),$$ as long as $q_i$
does not divide $A$. Thus, if $\Pi = 1-\frac {1}{A}$ cannot be
skipped by the averages of the sequences in
$SEQ(0,\mu_1,...,\mu_r,1)$ then $q_1|A,q_2|A,...,q_r|A$, that is,
$M|A$, or $$\Pi =1-\frac {1}{Mt}$$ for some $t\geq 1$ (the number
theory background necessary for the current paper can be found,
for example, in [2], Chapter 1). Thus we have proved that
$$IV(0,\mu_1,...,\mu _r, 1)\subseteq \left \{ 1-\frac {1}{Mt}:t=1,2,3,... \right
\}.\leqno (15)$$

\

\noindent To complete the proof we will prove the reverse
inclusion:

$$\left \{ 1-\frac {1}{Mt}:t=1,2,3,... \right
\} \subseteq IV(0,\mu_1,...,\mu _r, 1).\leqno (16)$$

\

\noindent We proceed by contradiction. Assume that $\frac
{Mt-1}{Mt}$ can be skipped by the consecutive averages of some
sequence in $SEQ(0,\mu_1,...,\mu_r,1)$. Without loss of generality
we may assume that $\frac {Mt-1}{Mt}$ is in between the average
$$\bar{x}_n= \frac {x_1+...+x_n}{n}$$
with $u$ of the $x_1,...,x_n$ being zeros, $v_1$ being $\frac
{p_1}{q_1}$, $v_2$ being $\frac {p_2}{q_2}$,..., $v_r$ being
$\frac {p_r}{q_r}$ and $w$ being ones ($u+v_1+...+v_r+w=n$), and
the average
$$\bar{x}_{n+1}= \frac {x_1+...+x_n+1}{n+1}$$
with $u$ of the $x_1,...,x_n,x_{n+1}=1$ being zeros, $v_1$ being
$\frac {p_1}{q_1}$, $v_2$ being $\frac {p_2}{q_2}$,..., $v_r$
being $\frac {p_r}{q_r}$ and $w+1$ being ones (we took $x_{n+1}=1$
which leads to the greatest possible increase in the average):
$$\frac {v_1\frac {p_1}{q_1}+...+v_r\frac {p_r}{q_r}+w}{n}<\frac {Mt-1}{Mt}
<\frac {v_1\frac {p_1}{q_1}+...+v_r\frac
{p_r}{q_r}+w+1}{n+1}.\leqno (17)$$ For every $i=1,...,r$ let us
define
$$Q_i:=\frac {{\rm lcm} (q_1,...q_r)}{q_i}=\frac {M}{q_i}.$$
With this notation, a multiplication of all terms in (17) by $M$
gives
$$\frac {p_1Q_1v_1+...+p_rQ_rv_r+Mw}{n} < \frac {Mt-1}{t} <
\frac {p_1Q_1v_1+...+p_rQ_rv_r+Mw+M}{n+1},$$ which, by
cross-multiplications turns out to be equivalent to the following
system of inequalities:
$$p_1Q_1v_1t+...+p_rQ_rv_rt+Mwt<Mnt-n,\leqno(18)$$
and
$$Mnt+Mt-n-1<p_1Q_1v_1t+...+p_rQ_rv_rt+Mwt+Mt\leqno (19)$$
Finally, from (18) and (19) it follows that
$$p_1Q_1v_1t+...+p_rQ_rv_rt+Mwt<Mnt-n<p_1Q_1v_1t+...+p_rQ_rv_rt+Mwt+1$$
which is, again, a contradiction (as there can be no integer
falling between consecutive integers). This shows that (16) is
true. From (15) and (16), (14) follows. This concludes the proof
of THEOREM 7.

\

\noindent From the above result, a linearity argument leads us to
the following arithmetic intermediate value theorem.

\

\

\noindent {\bf THEOREM 8.} Let $a_1<a_2<...<a_r$ ($r\geq 3$). For
$i=2,...,r-1$, let
$$\mu _i:=\frac {a_i-a_1}{a_r-a_1}.$$
Then the following hold true.

\

\noindent a) If for some $i=2,...,r-1$ the number $\mu _i$ is
irrational, then $IV(a_1,...,a_r)=\O$.

\

\noindent b) If $\mu_2,...,\mu _{r-1}$ are all rational numbers,
$\mu _i=\frac {p_i}{q_i}$, with $\gcd (p_i,q_i)=1$ for
$i=2,...,r-1$ and $M={\rm lcm}(q_2,...,q_{r-1})$, then
$$IV(a_1,...,a_r)=\left \{ \frac {1}{Mt}a_1+\left ( 1-\frac {1}{Mt}\right )a_r :t=1,2,3,...\right \}$$

\

\section {FURTHER COMMENTS}

\

\noindent Note that for $a_1<a_2<...<a_r$ the sets
$IV(a_1,...,a_r)$ represent the values $\Pi $ with the property
that the consecutive averages $\{ \bar{x}_n\}$ of every sequence
$\{ x_n\}_{n\geq 1}\in SEQ (a_1,...,a_r)$ cannot {\it increase}
from a value $\bar{x}_k<\Pi $ to a value $\bar{x}_l>\Pi $ without
taking the value $\bar{x}_s=\Pi$ for some $s$ with $k<s<l$.
Similarly we can define the sets $$DV(a_1,...,a_r)$$ representing
the the values $\Pi $ with the property that the consecutive
averages $\{ \bar{x}_n\}$ of every sequence $\{ x_n\}_{n\geq 1}\in
SEQ (a_1,...,a_r)$ cannot {\it decrease} from a value
$\bar{x}_k>\Pi $ to a value $\bar{x}_l<\Pi $ without taking the
value $\bar{x}_s=\Pi$ for some $s$ with $k<s<l$.

\

\noindent The connection between the sets $IV(a_1,...,a_r)$ and
$DV(a_1,...,a_r)$ can be expressed in a simple way as follows:
$$DV(a_1,...,a_r)=-IV(-a_r,...,-a_1).\leqno(20)$$
The proof of (20) is straightforward if we use the transformation
$$\{ x_n\}_{n\geq 1}\mapsto \{-x_n\}_{n\geq 1}.\leqno(21)$$
Clearly, (21) is a one-to-one correspondence between
$SEQ(a_1,...,a_r)$ and $SEQ(-a_r,...,-a_1)$. Under this
correspondence, the sequence of averages of $\{ x_n\}_{n\geq 1}$
skips (going up) $\Pi \in (a_1,a_r)$ if and only if the sequence
of averages of $\{-x_n\}_{n\geq 1}$ skips (going down) $-\Pi \in
(-a_r,-a_1)$.

\

\noindent We can use (21) to translate the Theorems 2 and 8 for
decreasing trends. Thus, we obtain

\

\noindent {\bf THEOREM 9.} If $a<b$, then
$$DV(a,b)=\left \{ \frac {k}{k+1}a+\frac {1}{k+1}b: k\geq 1  \right
\},$$

\

\noindent and

\

\noindent {\bf THEOREM 10.} Let $a_1<a_2<...<a_r$ ($r\geq 3$). For
$i=2,...,r-1$, let
$$\mu _i:=\frac {a_i-a_1}{a_r-a_1}.$$
Then the following hold true.

\

\noindent a) If for some $i=2,...,r-1$ the number $\mu _i$ is
irrational, then $DV(a_1,...,a_r)=\O$.

\

\noindent b) If $\mu_2,...,\mu _{r-1}$ are all rational numbers,
$\mu _i=\frac {p_i}{q_i}$, with $\gcd (p_i,q_i)=1$ for
$i=2,...,r-1$ and $M={\rm lcm}(q_2,...,q_{r-1})$, then
$$DV(a_1,...,a_r)=\left \{ \left ( 1-\frac {1}{Mt}\right )a_1+\frac {1}{Mt} a_r :t=1,2,3,...\right \}$$

\

\

\noindent {\bf REFERENCES}

\

\noindent [1] 2004 Putnam Exam, Problem A1

\noindent [2] I. Niven, H.S. Zuckerman and H.L. Montgomery, An
Introduction to the Theory of Numbers, $5$-th edition, Wiley 1991.

\

\end {document}